\documentclass[graybox]{SNmult}
\usepackage[utf8]{inputenc}
\usepackage[T1]{fontenc}

\usepackage{mathptmx}       % selects Times Roman as basic font
\usepackage{helvet}         % selects Helvetica as sans-serif font
\usepackage{courier}        % selects Courier as typewriter font
\usepackage{makeidx}         % allows index generation
\usepackage{graphicx}        % standard LaTeX graphics tool
\usepackage{multicol}        % used for the two-column index
\usepackage[bottom]{footmisc}% places footnotes at page bottom

\usepackage{newtxtext}       % 
\usepackage[varvw]{newtxmath}       % selects Times Roman as basic font

\usepackage[noadjust]{cite}
\usepackage{orcidlink}

\newtheorem{thm}{Theorem}%[section]
\newtheorem{cor}[thm]{Corollary}
\newtheorem{prop}[thm]{Proposition}

\begin{document}

\title*{Topological aspects of generalized Sklyanin algebras}
\titlerunning{Topological aspects of generalized Sklyanin algebras}

\def\oorcid{0000-0000-0000-1111}

\author{Giorgi Khimshiashvili\orcidID{0000-0001-8056-1120}}
\authorrunning{Giorgi Khimshiashvili}

\institute{Giorgi Khimshiashvili\at
  Faculty of Business,Technology and Education\\
  Ilia State University\\
Ilia Chavchavadze Ave. N45\\
Tbilisi 0179, Georgia\\
\email{giorgi.khimshiashvili@iliauni.edu.ge}\orcidlink{https://orcid.org/0000-0001-8056-1120}}

\maketitle

\abstract{We are concerned with the topology of Casimir foliations of
  Sklyanin algebras and their certain generalizations. In particular,
  we give a complete description of the topological types and
  bifurcations of Casimir leaves of the original Sklyanin algebras,
  which improves and extends analogous results from the original
  papers of E.~Sklyanin. Our approach is based on the so-called
  signature formula for the Euler characteristic of real algebraic
  varieties. Moreover, we show that the same approach yields similar
  results for certain generalizations of the original Sklyanin
  algebras known as Jacobian Poisson structures.}

\keywords{Poisson structure, Sklyanin algebra, Casimir function,
  Casimir foliation, Nambu bracket, Jacobian Poisson structure, Euler
  characteristic, algebraic formula for the Euler characteristic,
  bifurcation diagram, isolated singularity of complete intersection,
  Milnor number\vspace*{5pt}}


\section{Introduction}\label{S1}

We deal with the so-called Sklyanin algebras introduced in the seminal papers of E.~Sklyanin \cite{skl1}, \cite{skl2} and their generalizations studied in \cite{khi1}, \cite{khpr}, \cite{prz1}. As explained in \cite{skl2} the topology of the so-called {\it Casimir foliation} of Sklyanin algebra is closely related to the representations of the quantized Sklyanin algebra, which was the central topic in \cite{skl2}. However the results on the topology of Casimir leaves in \cite{skl1}, \cite{skl2} were presented without detailed proofs. Moreover, the description of the topology of Casimir levels given in \cite{skl1} contained a slight error which was corrected in \cite{skl2} but the proof of the correct result was only sketched. Later on, the present author realized that the topological results given in \cite{skl1}, \cite{skl2} can be rigorously proved and generalized using standard tools of algebraic geometry and nonlinear analysis. The aim of the present paper is to present some of the results obtained in this way including the rigorous proofs of the mentioned results of E.~Sklyanin. 

\section{Sklyanin algebras and Jacobian Poisson brackets}\label{sec:2}

We begin by presenting the necessary background on Sklyanin algebras and their generalizations known as Jacobian Poisson structures. The Sklyanin algebras were introduced in \cite{skl1} as Poisson brackets on the algebra of polynomial functions $\mathbb{R}_4$ on four-dimensional Euclidean space $\mathbb{R}^4$ defined by the following relations 

\begin{equation}\label{sklyanin1} 
	\{x_0, x_i\} = J_{jk}x_jx_k, \{x_i, x_j\}= -2x_0x_k. 
\end{equation}
Here $x_0, x_1, x_2, x_2$ are the coordinate functions, $J_1, J_2, J_3$ are non-zero real numbers, and $J_{jk}=J_j-J_k$ for a cyclic permutation $(i,j,k)$ of $(1,2,3)$. For convenience and to comply with the original definition it will be assumed that $J_1 \geq J_2 \geq J_3$. We preserve the original notation of E.~Sklyanin to facilitate the comparison of our results with those given in \cite{skl1}, \cite{skl2}.

As was shown in \cite{skl1}, these Poisson brackets possess two quadratic Casimir polynomials $K_0, K_1$ given by   

\begin{equation}\label{casimirs1} 
	K_0 = x_1^2 + x_2^2 + x_3^2, K_1 = x_0^2 + J_1x_1^2 + J_2x_2^2 + J_3x_3^2.
\end{equation}

The quantized versions of these algebras are related to various exactly solvable problems of quantum mechanics.
As was explained in \cite{skl2}, certain information about the number and types of their representations could be obtained by studying the topology and bifurcations of the so-called {\it Casimir leaves} defined as the level surfaces $X_{ab} = \{K_0=a, K_1=b\}$ of Casimir functions $K_0, K_1$. For this reason, several results on the topology of Casimir leaves were given in \cite{skl1} but their proofs were only outlined. Later on, it turned out that the description of the topological structure of Casimir leaves presented on page 30 of \cite{skl1} contains an error. This error was corrected in \cite{skl2} but again without a detailed proof. Moreover, the description of bifurcations in the Casimir foliation remained incomplete. One of the aims of this paper is to give rigorous proofs of the mentioned results and give a more detailed description of bifurcations in the Casimir foliation of a Sklyanin algebra.

For the sake of clarity, we outline now our approach and the structure
of the paper. We consider the mapping
$K: \mathbb{R}^4 \rightarrow \mathbb{R}^2$ defined by the Casimir
polynomials $K_0, K_1$, show its properness and study its fibers using
methods of nonlinear analysis and real algebraic
geometry. Specifically, we compute bifurcation diagram $B(K)$ and
compute the Euler characteristic of non-singular fibers over each {\it
  chamber} (connected component of the complement of $B(F)$ in
$\mathbb{R}^2$) using the so-called signature formula for the Euler
characteristic developed in our earlier works (see, e.g., \cite{khi1},
\cite{khi2}). A computer-assisted computation showed that the Euler
characteristic of non-singular Casimir leave can only be equal to zero
or four. Since non-singular fibers are compact two-dimensional
surfaces it follows that they are diffeomorphic either to a
two-dimensional torus $T^2$ or to a disjoint union of two
two-dimensional spheres $S^2$. This, in particular, shows that there
are no fibers diffeomorphic to sphere $S^2$ as was erroneously claimed
on the page~30 of~\cite{skl1}.

Then we compute in the same way the Euler characteristic of fibers over each component of $B(K)$ and find out that it is always equal to two. In order to study bifurcations of Casimir leaves we invoke a well known fact that quadratic mapping into the plane has only singularities of quadratic cone type \cite{agr}. This implies that all singular fibers are homeomorphic to a torus $T^2$ with two pinched meridians, which can be informally described as the union of two bananas with identified tips. These results are summed up in Theorem \ref{main1} below. In the rest of the paper we present similar results for generalized Sklyanin algebras defined by Jacobian Poisson brackets (Theorem \ref{main2}) and outline further possible developments in the same spirit. 

We conclude this section by recalling the definition of generalized Sklyanin algebras in terms of Jacobian Poisson structures which play the main role in Section 4. To this end we recall first the definition of {\it Jacobian Poisson structures} which were introduced in \cite{khi3} and studied in some detail \cite{khpr} and \cite{prz1}. 
Given a collection $F$ of $n-2$ functionally independent real polynomials $F=\{f_1, \ldots, f_{n-2}\}$ of $n$ variables the Jacobian Poisson brackets (JPB) of two smooth functions $g, h$ on $\mathbb{R}^n$ are defined by the formula

\begin{equation} \label{Jbrackets}
	\{g,h\}_F = J(f,g,f_1, \ldots, f_{n-2}),
\end{equation}
where $J$ stands for the jacobian of the functions $f,g,f_1, \ldots, f_{n-2}$. Functions $f_j$ are called the {\it generators} of JPB (\ref{Jbrackets}). Formula (\ref{Jbrackets}) is a particular case of the so-called {\it Nambu brackets}, which implies that $\{g,h\}_F$ defines a genuine Poisson structure on $\mathbb{R}^n$ called a {\it Jacobian Poisson structure} (JPS). The basic properties of JPS have been established in \cite{khi3} and \cite{khpr}.
The generators $\{f_1, \ldots, f_{n-2}\}$ are obviously Casimir functions of JPS (\ref{Jbrackets}), which define a foliation called the Casimir foliation of JPS (\ref{Jbrackets}). The non-singular leaves of this foliation are two-dimensional as in the case of Sklyanin algebras. The following proposition from \cite{khi3} shows that Sklyanin algebras are a special case of JPS and suggests several natural generalizations of Sklyanin algebras some of which are discussed in the sequel.

\begin{prop} \label{sklyanin-nambu1}
	Sklyanin algebra is a JPS defined by Jacobian Poisson brackets (\ref{Jbrackets}) with generators $K_0, K_1$:
	\begin{equation} \label{nambu1}
		\{f, g\}_S = J(f, g, K_0, K_1).
	\end{equation}
\end{prop}

If all generators of JPS  (\ref{Jbrackets}) are homogeneous quadratic polynomials (quadratic forms) and all leaves of its Casimir foliation are compact then such a JPS will be called a {\it generalized Sklyanin algebra} (GSA). There also exist other generalizations of Sklyanin algebras but in the sequel this term is only used in the above sense. 
In the rest of the paper we are basically concerned with the topology of Casimir foliations of Sklyanin algebras and 
generalized Sklyanin algebras. A few of other generalizations will be briefly discussed from the same perspective in the last section of the paper.  

From the defintion of GSA given above follows that its non-singular Casimir leaves are smooth compact two-dimensional surfaces. Hence important information on the topology of Casimir leaves and Casimir foliation is given by the values of Euler characteristic of its fibers. For this reason, the main tools in our approach are algebraic formulas for the Euler characteristic given in the next section.

\section{Algebraic formulas for the Euler characteristic}
\label{sec:3}

We now present the formulas for the Euler characteristic used in the sequel. The first one is an algebraic formula for the Euler characteristic of the so-called {\it real Milnor fibre} $\{f=0\}\cap S^{n-1}_r$ of an isolated singularity of real algebraic hypersurface $\{f=0\}$ in $\mathbb{R}^n$, where $S^{n-1}_r$ is a hypersphere of small radius $r$ centered at the origin. It is known that, for sufficiently small $r>)$ this intersection is a smooth algebraic surface of dimension $n-2$. It follows that, for an odd $n$, the Euler characteristics of this set vanishes while, for an even $n$, we have the following formula which was proved in \cite{khi1}. 

\begin{thm} [\cite{khi1}]\label{chi1}
	For even $n$,
	\begin{equation} \label{chiformula1}
		\chi(\{f=0\}\cap S^{n-1}_r) = 2 - 2\operatorname{deg}_0 \nabla f,
	\end{equation}
	where $\operatorname{deg}_0 \nabla f$ is the local topological degree of the gradient $\nabla f$ (aka the index of the gradient vector field).
\end{thm} 

The above formula will enable us to compute the Euler characteristics of Casimir leaves of Sklyanin algebra. In the case of generalized Sklyanin algebras we need a more general formula for the Euler characteristic of compact real algebraic set $X$ which was proven in \cite{khi2}, \cite{bru}. Suppose that compact real algebraic set $X$ is given as the set of common zeroes of $k$ polynomials $f_j(x_1, \ldots, x_n)$ of degree $d$ in $\mathbb{R}^n$. Let us introduce the auxiliary polynomials $h_j$ $\mathbb{R}^{n+1}$ by the formulas 
$$h_j (x_0, x_1, \ldots, x_n) = x_0^{d+1} f_j(\frac{x_1}{x_0}, \ldots, \frac{x_n}{x_0}), j=1, \ldots, k.$$ After that we define the "master polynomial" $H$ by

	\begin{equation} \label{Hpolynomial}
	H = \sum_{j=1}^n h_j^2 - \sum_{j=0}^n x_j^{2d+4}.
\end{equation}

\begin{thm} [\cite{khi2}, \cite{bru}] \label{chi2} With the above
  assumptions and notations, polynomial $H$ has an isolated critical
  point at the origin and
  \begin{equation} \label{chiformula2} 2 \chi(X) = 1 - \operatorname{deg}_0 \nabla H,
  \end{equation}
  where $\operatorname{deg}_0 \nabla H$ is the local topological degree of the
  gradient $\nabla H$.
\end{thm}

The above formulas being conceptually important are moreover practically applicable due to the so-called signature formula for the topological degree established in \cite{khi1}, \cite{eile}. The latter formula yields an algorithmic way of computing the local gradient degree participating in formulas (\ref{chi1}, \ref{chi2}). In the next section, we use formulas (\ref{chi1}, \ref{chi2}) to compute the Euler characteristics of Casimir leaves of Sklyanin algebras.
To do that in a rigorous way we need to ensure certain topological properties of the Casimir mapping, which is done in the next section.

\section{Casimir foliation of Sklyanin algebra}
\label{sec:4}

We give now the necessary topological properties of Casimir mappings and outlines of their proofs.
These properties are basically corollaries of the following three properties of quadratic mappings
into the plane which are well known (see, e.g., \cite{khi2}). 
\begin{description}
\item [\textbf{P1.}] Quadratic mapping is proper if and only if the
  pencil of quadrics generated by its components contains a positively
  determined quadratic form.
\item [\textbf{P2.}] The image of quadratic mapping into the plane is
  convex.
\end{description}

\begin{prop} \label{casimirproper}
	For any Sklyanin algebra, its Casimir mapping is proper. 
\end{prop}

This follows from P1 by noticing that the sum $MK_0 + K_1$ is a positively determined quadratic form for sufficiently
big $M>0$. 

\begin{prop} \label{casimirconvex}
	For any Sklyanin algebra, the image ot its Casimir map is convex.
\end{prop}

This is a direct corollary of P2. Thus the above two properties are special cases of general results on quadratic 
mappings into the plane. The next two properties follow from the explicit formulas
of Casimir polynomials of Sklyanin algebra.

\begin{prop} \label{casimirimage}
	For any Sklyanin algebra, the image ot its Casimir mapping is a halfplane. 
\end{prop}

Indeed, it is easy to verify that, for any pair of numbers $(u, v)$ with $u>0$,  the system $K_0 = u, K_1 = v$ has real solutions.  

\begin{prop} \label{casimirsdiscrim}
	For any Sklyanin algebra, the discriminant of its Casimir mapping consists of three rays in the upper half-plane.   
\end{prop}

This follows by explicit computation of the singular set $S(K)$ of Casimir mapping and its discriminant $B(K)= K(S(K))$. More precisely, using the notations of Section~2 it can be verified that $B(K)$ consists of three rays 
$v = J_ju, j =1, 2, 3,$ in the upper half-plane $u\geq 0$. 

This done we can describe the topology of Casimir foliation as follows. Let us denote by $U_j$ the four sectors which
are the components of the complement of $B(K)$ in the half-plane $u\geq 0$. By Ehresmann theorem the topological type of the preimage $K^{-1}(P)$ does not change within each $U_j$. Taking a point $P_j$ in each $U_j$ and computing 
the Euler characteristic $\chi(K^{-1}(P_j))$ by formula (\ref{chi2}) one finds that it only takes two values - zero and four. Since the Casimir leaves are smooth compact two-dimensional surfaces it follows that each Casimir leave is diffeomorphic either to the torus $T^2$ or to a disjoint union of two two-dimensional spheres $S^2$. The same computation shows that the leaves of the second type correspond to the outer sectors adjacent to the axis $u=0$, while the toric leaves correspond to the two inner sectors. It remains to describe the topological type of the leaves over each of the rays $u = J_jv$. This is again possible using formula (\ref{chi2}), which yields that the Euler characteristic of singular leaves is equal to two. Recalling that the singular points of the fibers of proper quadratic mappings into the plane are od quadratic cone type we conclude that the bifurcation happens by pinching two meridians of the toric leave as was claimed in the introduction. One can now describe the Reeb graph of the Casimir mapping in the spirit of \cite{van} and give a more detailed picture of bifurcations by explicit formulas. 
The results obtained are summed up in the following theorem.

\begin{thm} \label{main1}
	The complete list of homeomorphic types of Casimir leaves is: disjoint union of two two-dimensional spheres, two-dimensional torus, and torus with two meridians pinched. There is only type of bifurcation in Casimir foliation which corresponds to pinching two meridians of a two-torus. 
\end{thm}

In the next section we present a partial generalization of Theorem (\ref{main1}) to the case of GSA.  We add that the same results can be obtained using the so-called {\it Agrachev spectral sequence} (see, e.g., \cite{agr}). 

\section{Topology of generalized Sklyanin algebras}
\label{sec:5}

The discussion in this section follows the same way of reasoning. We establish first a characterization of GSA in terms of properties of its generators. 

\begin{thm}\label{main2}
	A collection $Q$ of $n-2$ real $n$-ary quadratic forms $q_j$ can serve as generators of a GSA $\mathbb{R}^n$ if and only if their Jacobi determinant has rank $n-2$ almost everywhere, and there exists a positively definite linear combination of the forms $q_j$. The image of Casimir mapping $K_Q: \mathbb{R}^n \rightarrow \mathbb{R}^2$ is a halfplane and its discriminant $B(K_Q)$ is a collection of smooth curves emanating from the origin of $\mathbb{R}^2$.
\end{thm}

This follows from the above discussion in the same way as Propositions (\ref{casimirproper}), (\ref{casimirconvex}), (\ref{casimirimage}) combined with a direct computation of the discriminant. In this setting the number of  discriminant branches depends on the generators $q_j$ and we are only able to obtain an upper estimate of this number. The emerging estimate of Bezout type is far from being exact and therefore we omit it. An upper estimate for the number of possible topological types of Casimir leaves can be given in terms of the so-called {\it Petrovsky numbers} introduced by V.Arnol'd in \cite{arn1}. 

Recall that Petrovsky number $P(M, m)$ was defined in \cite{arn1} as the number of $n$-ary monomials of degree $m(M-2)$. According to \cite{arn1} the number of possible values of the gradient degree of polynomial of degree $M$ in $m$ variables does not exceed $2P(M, m)$. Moreover, formula (\ref{chi2}) implies that the number of possible values of the Euler characteristic of Casimir leaves of GSA does not exceed the number of possible values of the gradient degree of the master polynomial $H$ given by formula (\ref{Hpolynomial}). In case of GSA, the master polynomial depends on $n+1$ variables and has degree $8$. Hence the relevant Petrovsky number is $P(8, n+1)$, 
which implies that the number of non-diffeomorphic Casimir leaves does not exceed
$$\frac{(n+6)!}{8!(n-2)!}.$$  

This estimate is obviously far from being exact but it can be improved if the generators $q_j$ are simultaneously diagonalizable in a certain coordinate system as happens in the case of original Sklyanin algebra. An interesting problem is to improve the above estimate for various special classes of GSA but this is beyond the aims of the present paper.   

We proceed by noting that in addition to the above direct generalizations of the results on Sklyanin algebras, one can consider several further topological aspects of GSA involving more refined concepts and technical tools. As an example, we notice that the definition of GSA given above implies that the collection of its quadratic Casimirs defines a so-called {\it isolated singularity of complete intersection} $X = \{F=0\}$ (ICIS) in the sense of singularity theory. One of the most important topological invariants of such a singularity is its Milnor number (see, e.g., \cite{loo}). Using it we can obtain certain information on the topology of Poisson structure induced by GSA on the common zero-set $X$ of its Casimirs.

\begin{thm} \label{milnorGSA}
	For any generalized Sklyanin algebra on $\mathbb{R}^n$ with $n\geq 4$, the dimension of zero-th Poisson cohomology of the induced Poisson structure on the zero-set of Casimirs is equal to
	
	\begin{equation} \label{milnorGSA}
		\sum_{j=0}^{n-3} (-1)^{n-j-3}\frac{n!}{j!(n-j)!}2^j. 
	\end{equation}  
\end{thm}

This follows from the following two general results. First, as was proved in \cite{khi3}, for any JPS, the dimension of zero-th Poisson cohomology of the induced Poisson structure on the zero-set of its Casimirs is equal to the Milnor number of the ICIS defined by the generators. The second general fact is the following formula for the Milnor number of an ICIS defined by homogeneous polynomials of arbitrary degrees, which was given by G.Greuel and H.Hamm in \cite{grha} (see  \cite{loo} for a detailed and simplified proof).

\begin{thm} (\cite{grha}) \label{milnorICIS}
	\begin{equation} \label{milnorICIS}
		\mu(f) = \sum_{j=0}^{k-1} (-1)^{k-j-1}\frac{n!}{j!(n-j)!}2^j. 
	\end{equation}
\end{thm}  

In view of the said above, formula (\ref{milnorGSA}) is obtained by merely combining the above observations with formula (\ref{milnorICIS}) and taking into account that $k=n-2$ in the case of a generalized Sklyanin algebra on $\mathbb{R}^n$. 

\begin{cor}
	For an original Sklyanin algebra, the dimension of zero-th Poisson cohomology of the induced Poisson structure on the zero-set of Casimirs is equal to seven. 
\end{cor}

It would be interesting to give explicit formulas for the generators of the zero-th Poisson cohomology of Sklyanin algebras and generalized Sklyanin algebras.  

\section{Concluding remarks}
\label{sec:6}

The approach and results presented above suggest several further
developments along the same lines. For example, the upper estimate for
the number of possible topological types of the Casimir leaves of GSA
in terms of Petrovsky number can be improved for special classes of
GSA. In particular, if $n=5$ and all Casimir polynomials admit
simultaneous diagonalization, all possible topological types have been
listed by V.~Krasnov \cite{kra}.

Next, on more sophisticated level, an obvious research prospect is to
compute the Poisson (co)homology of generalized Sklyanin algebras. As
a paradigm one can use the comprehensive results on the Poisson
cohomology of the original Sklyanin algebras given by S.R.T.~Pelap
\cite{pel}. A natural next step would be to extend the results of
\cite{pel} to all generalized Sklyanin algebras in dimension four.

Finally, analogs of Theorem \ref{milnorGSA} can be obtained for Jacobian Poisson structures with Casimirs of arbitrary degree using the general formula (\ref{milnorICIS}).  

Developments suggested above require separate discussion and will (hopefully) be considered elsewhere.

\end{document}